\documentclass[pdflatex]{sn-jnl}

\usepackage{graphicx}%
\usepackage{multirow}%
\usepackage{amsmath,amssymb,amsfonts}%
\usepackage{amsthm}%
\usepackage{mathrsfs}%
\usepackage[title]{appendix}%
\usepackage{xcolor}%
\usepackage{textcomp}%
\usepackage{manyfoot}%
\usepackage{booktabs}%
\usepackage{algorithm}%
\usepackage{algorithmicx}%
\usepackage{algpseudocode}%
\usepackage{listings}%
\usepackage{tikz}%
\usetikzlibrary{arrows.meta}%

\theoremstyle{thmstyleone}%
\newtheorem{theorem}{Theorem}%

\theoremstyle{thmstyletwo}%

\theoremstyle{thmstylethree}%
\newtheorem{definition}{Definition}%

\begin{document}

\title[Three Wrong Definitions]{The History of Three Wrong Definitions}

\author{\fnm{Harold P.} \sur{Boas}}\email{boas@tamu.edu}

\affil{\orgdiv{Department of Mathematics}, \orgname{Texas A\&M University}, \orgaddress{\city{College Station}, \postcode{77843-3368}, \state{TX}, \country{USA}}}

\abstract{The topic is the history of the concepts of equivalence relation, Cauchy sequence, and metric space. The thesis is that disused definitions of these notions could profitably be revived.}

\keywords{equivalence relation, Cauchy sequence, metric space}

\pacs[MSC Classification]{01-01, 03-03, 40-03, 54-03}

\maketitle

\vspace*{-24pt}
\hrulefill

\noindent 
This document is the author's Accepted Manuscript, accepted 3~March 2026 after peer review for publication in \emph{The Mathematical Intelligencer}. It is
not the Version of Record and does not reflect post-acceptance improvements or any corrections. The Version of Record, published online 8~April 2026, is
available at
\url{https://doi.org/10.1007/s00283-026-10520-7}.

\noindent
\hrulefill

\bigskip\bigskip

I need a volunteer to help perform a magic trick. Are you game?

Your task is to pose a certain 
mathematical question to a third party and to record the answer (without letting me see).
If you are a student, you might text the question to a classmate. If you are an instructor, you might poll a class or consult a colleague. Alternatively, you may query your favorite artificial-intelligence assistant. 

The question is based on the background of your interlocutor. 
For a beginner just learning how to write proofs, the question is to specify the defining properties of an equivalence relation. For someone familiar with rigorous calculus, the question is to define the concept of a Cauchy sequence of real numbers. If your informant has studied the elements of point-set topology, the question is to state the axioms for a metric space. 

After obtaining an answer, read on for my divination of what you were told. Imagine that I have made a magic incantation ending with ``Eureka!''

If you asked the first question, the response you got was that an equivalence relation satisfies the following three properties (implicitly universally quantified).
\begin{itemize}
    \item Reflexivity: \(x\) is related to~\(x\).
    \item Symmetry: if \(x\) is related to~\(y\), then \(y\)~is related to~\(x\).
    \item Transitivity: if \(x\) is related to~\(y\), and \(y\)~is related to~\(z\), then \(x\)~is related to~\(z\).
\end{itemize}

If you asked the second question, the reply was that a sequence \(\{x_n\}_{n=1}^\infty\) of real numbers is a Cauchy sequence when
\begin{equation}
(\forall \varepsilon>0) (\exists N) (\forall m\ge N) (\forall n\ge N) ( |x_m-x_n|<\varepsilon).
\label{eq:Cauchy}
\end{equation}

The answer to the third question was that a metric space means a set~\(E\) together with a nonnegative function \(d\) on \(E\times E\) satisfying the following three properties, again implicitly universally quantified. 
\begin{align}
&d(x,y)=0 \quad \text{if and only if} \quad x=y \quad \text{(coincidence axiom)}.
\label{eq:id} \\
&d(x,y)=d(y,x) \quad \text{(symmetry axiom)}. \label{eq:sym} \\
&d(x,z) \le d(x,y) + d(y,z) \quad \text{(triangle inequality)}.
\label{eq:tri}
\end{align}

Did my mind-reading trick succeed?
The magic here was solely in the selection of the questions. If the topics had been instead function, real number, and vector space, then I could not have confidently predicted the response. 
The three concepts that I chose, however, were codified in textbooks  long ago, and modern students all learn the definitions in substantially the same conventional form.  
The statements are a part of mathematical culture, like
the legend of Archimedes
in the bath.

The tale of Archimedes contains a grain of truth but is apocryphal. 
My startling assertion is that analogously, these three standard textbook definitions are plausible but wrong! 

To be sure, I am not suggesting 
that these definitions are unsound. 
What I mean is that the ``right'' definitions are shorter, simpler, and more insightful. 
In audaciously challenging the conventional wisdom,
I am being playfully provocative, in the spirit of ``\(\pi\) is wrong!''\footnote{A quarter century ago, Bob Palais promoted in this journal \cite{palais} the opinion that the right circle constant is the ratio of the circumference to the radius, namely, \(2\pi\) rather than~\(\pi\).} 

My initial goal in this article is to formulate and to popularize some appealing alternatives to the prevailing definitions of the three concepts. Dating back a century or more, these formulations curiously have dropped out of the everyday toolbox of mathematical workers.

My second objective is to trace the history in sufficient detail to name the creators. Today's standard statements seemingly arose as much by chance as by choice, and
I try to identify the historical turning points when the variant definitions slipped through the fingers of Fame and fell 
into footnotes. 

My third aim is to light candles in remembrance of scholars whose names, like their definitions, are fading from collective consciousness. 
Both for obscure figures and for revered icons, I offer glimpses into their lives to enhance the appreciation of their achievements. 
For example, 
Augustus De~Morgan's recipe for the algebra of sets is more piquant when seasoned with the awareness that as a boy blind in one eye, De~Morgan overcame bullying at school.

\section*{All things being equal}
Did the ancient Greeks anticipate the idea of an equivalence relation? 
If Euclid is accepted as a representative of classical mathematics, then the evidence is ambiguous but tantalizing. 

The first ``common notion'' in Euclid's \emph{Elements} states,
``Things which are equal to the same thing are also equal to one another'' (Heath's translation~\cite[p.~155]{heath}). 
Although Euclid never defines equality \cite{lucic}, 
the import evidently is more general 
than simple identity. If ``equal'' is reinterpreted as ``equivalent,'' then Euclid's first common notion becomes the property used to prove that two equivalence classes must be either identical or disjoint. 

The \emph{Elements} calls congruent line segments ``equal.'' Robin Hartshorne observes, ``there are no lengths in Euclid's geometry, so \ldots\ we may regard equality \ldots\ to be an equivalence relation on line segments'' \cite[p.~28]{hartshorne}.
 
Book~I of the \emph{Elements} culminates in the Pythagorean theorem, which states that the square on the hypotenuse of a right triangle ``is equal to'' the squares on the other two sides. Euclid's meaning is formally stronger than equality of areas, for Euclid calls rectilinear figures ``equal'' if they can be converted into congruent figures through the operation of removing or adding congruent triangles (iterated finitely many times). 
This notion of equality is an equivalence relation.

Book~V of the \emph{Elements}, presumed to be based on lost work of Eudoxus of Cnidus, defines ratios \(a:b\) and \(c:d\) of abstract quantities to be ``the same'' when for arbitrary natural numbers \(m\) and~\(n\), the inequality \(ma > nb\) holds if and only if \(mc > nd\), and similarly \(ma<nb\) if and only if \(mc < nd\). 
Does the device of declaring distinct objects to be ``the same'' ring a bell? In the words of Christopher Zeeman (1925--2016), ``Eudoxus and Euclid must have been thinking of a ratio as something like an equivalence class'' \cite[p.~16]{zeeman}.

Thus Euclid knew some particular examples of mathematical structures that now are understood as equivalence relations. According to
the principle of indeterminacy of translation \cite[chap.~2]{quine} of W.~V.~O. 
Quine (1908--2000), there is no way to know how the ancients themselves conceptualized these examples. I doubt that Euclid had any notion of abstract equivalence relations in general, for the relation concept did not get formalized by mathematicians for another two millennia. 

Adrien-Marie Legendre (1752--1833) published a textbook in 1794 intended to make Euclid's geometry more accessible to students of the time.\footnote{Legendre's influential book went through many editions during his lifetime and multiple translations and posthumous revisions by various editors. Yet Lewis Carroll (1832--1898) maintained, based on extensive teaching experience, that the book was suitable for advanced students only, not for beginners \cite[p.~59]{carroll}.} 
Reserving the word ``equal'' for planar figures that are congruent,  he introduced the word ``equivalent'' for figures that have equal areas \cite[p.~57]{legendre}. During the ensuing 19th century, the concept of an equivalence class was ``in the air'' \cite[p.~4663]{asghari}.

This atmospheric disturbance produced multiple mathematical lightning strikes. The first bolt from the blue was the publication in 1801 of a seminal treatise on number theory in which the young Carl Friedrich Gauss (1777--1855) developed the basic example of congruence of integers~\cite{gauss}. 
The thunderclap still reverberated half a century later in the ears of Augustin-Louis Cauchy (1789--1857), who proposed an innovative definition of the complex numbers~\(\mathbb{C}\) in terms of equivalence classes of polynomials, \(\mathbb{R}[x]/(x^2+1)\) in modern notation~\cite{cauchy-quotient}. 
In 1871, 
Richard Dedekind \cite{sonar} (1831--1916) 
introduced an equivalence relation on ideals in the ring of integers of an algebraic number field and 
initiated the study of the ideal class group
\cite[\S~164]{dirichlet}. In the same decade, 
Georg Cantor (1845--1918) 
defined infinite cardinalities as equivalence classes under the relation of bijection~\cite{cantor}. 

Although these Prometheans revealed the secrets of specific equivalence relations, they left the formulation of a general definition to mere mortals.
Only in the second half of the 19th century did a theory of abstract binary relations gradually precipitate out of a fog of uncertainty surrounding the foundations of mathematical logic. 

In 1850, Augustus De Morgan (1806--1871) highlighted the importance of the pair of properties that he called ``transitiveness'' and ``convertibility'' \cite[p.~104]{DeMorgan} (now ``transitivity'' and ``symmetry''). 
A relation satisfying these two special properties is an ``abstract copula'' in his terminology.
Charles Sanders Peirce (1839--1914) used block-diagonal matrices a quarter century later to describe ``copulatives'' \cite[p.~368]{peirce},
which are recognizable as being the equivalence relations. Since 
Peirce\footnote{The most original thinker in an illustrious family, C. S. Peirce should not be confused with relatives who held posts at Harvard. His elder brother, James Mills Peirce (1834--1906), was the second Perkins Professor of Mathematics. These brothers' father was the renowned Benjamin Peirce (1809--1880), the first holder of the Perkins chair and the eponym of present postdoctoral fellowships. A third cousin of the brothers was Benjamin Osgood Peirce (1854--1914), the seventh Hollis Professor of Mathematicks and Natural Philosophy. (Quine was the Edgar Pierce Professor of Philosophy, but this unrelated Pierce is spelled differently.)} 
was working in America, then a mathematical backwater, his ideas were slow to influence the current of mathematical thought in Europe.

Indeed, a calculus textbook of 1888 by Giuseppe Peano (1858--1932) says that a symmetric and transitive relation is an ``equality''
\cite[p.~141]{calcolo}, wording reminiscent of Euclid. 
When Peano formulated his axioms for the natural numbers the following year, he supplemented the definition of an equality with the assumption that every element is related to itself \cite[p.~XIII]{arithmetices}. Peano's disciple Giovanni Vailati (1863--1909) suggested naming this additional property ``reflexivity'' 
\cite[p.~134]{vailati-reflex}, a designation that was universally adopted.

The reflexive, symmetric, transitive relations might well have ended up being called ``equality relations,'' for Peano's terminology was still  circulating as late as 1929, when the philosopher Rudolf Carnap (1891--1970) published his survey of symbolic logic \cite[p.~48]{carnap}. The neologisms ``equiparative relation'' \cite[p.~113]{deamicis}, due to 
Enrico De Amicis (1858--1925),
and ``isoid relation'' \cite{jourdain}, due to Philip E.~B. Jourdain (1879--1919), gained no traction. 

Helmut Hasse (1898--1979) gets credit for the name \emph{\"Aquivalenzrelation}. In addition to making fundamental contributions in algebra and number theory, he succeeded his teacher Kurt Hensel (1861--1941) as editor of the venerable Crelle journal (\emph{Journal f\"ur die reine und angewandte Mathematik}), a post in which Hasse served for half a century~\cite{rohrbach}. 
Hasse standardized the following definition in the first volume of his treatise on algebra \cite[I.2.2]{hasse}, published one hundred years ago.
\begin{definition}
\label{def:hasse}
    An equivalence relation is a relation that is reflexive, symmetric, and transitive. 
\end{definition}

Hasse's definition sounds perfectly in tune, as to be expected from the hands of a skilled pianist. The three indicated properties are independent of each other, as Peano had observed decades earlier \cite[p.~24]{geometria}, so the definition apparently forms a harmonious chord. Nonetheless, two discordant overtones can be detected. 

First, the character of reflexivity contrasts with the other properties. The reason that reflexivity failed to catch De~Morgan's attention is, I imagine, that symmetry and transitivity are \emph{implications}, but reflexivity is a mandate.
Moreover, the following discussion shows that reflexivity is almost automatic in the presence of symmetry and transitivity. 

In the terminology introduced by Bertrand Russell\footnote{Russell was pictured recently in this journal in the Stamp Corner \cite{bisht}.} 
(1872--1970) at the turn of the 20th century \cite[p.~116]{russell}, an element~\(x\) belongs to the \emph{domain}\footnote{Russell's parallel term \emph{converse domain} for the domain of the converse (or inverse) relation was promptly compressed by his follower Louis Couturat (1868--1914) to \emph{codomain} \cite[p.~29]{couturat}.\label{foot:couturat}} of a relation~\(\sim\) on a set when there exists an element~\(y\) for which \(x\sim y\). If the relation~\(\sim\) is symmetric, then also \(y\sim x\). When \(\sim\) additionally is transitive, the properties that \(x\sim y\) and \(y\sim x\) imply that \(x\sim x\). This remark (which goes back to the 19th century~\cite[p.~127]{deamicis}) yields the first variation on Hasse's definition.
The essence of the following statement can be found as a proposition in Bourbaki's treatise on set theory \cite[p.~II.41]{bourbaki}. 

\begin{definition}
An equivalence relation on a set is a symmetric and transitive relation whose domain is the whole set. 
\end{definition}

The second dissonance in Definition~\ref{def:hasse} arises because symmetry deals with pairs of elements, and transitivity deals with triples. Each pair arises from some triple by setting two elements equal, so every property of pairs is a special case of some property of triples. Is there a natural property of triples that specializes to symmetry?

\begin{figure}
    \centering
 \begin{tikzpicture}[>=Stealth,thick]
    \fill (0,0) node[anchor=north east]{$x$} circle (2pt);
    \fill (2,3.464) node[anchor=south] {$y$} circle (2pt);
    \fill (4,0) node[anchor=north west] {$z$} circle (2pt);
  \draw[->] (0,0) +(60:2pt) -- (60:4cm-2pt) ;    
    \draw[->] (60:4cm) +(300:2pt) -- +(300:4cm-2pt);   
    \draw[->,dashed] (2pt,0pt) --(4cm-2pt,0pt); 
 \end{tikzpicture}
    \caption{Transitive relation}
    \label{fig:transitive}
\end{figure}
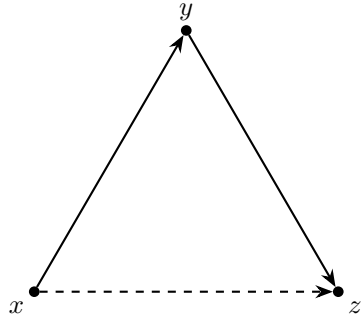

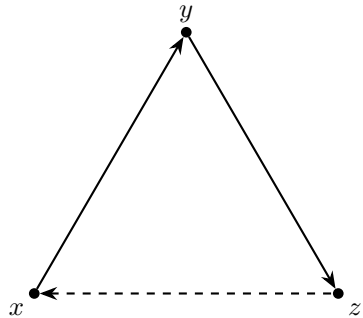
\begin{figure}
    \centering
 \begin{tikzpicture}[>=Stealth,thick]
    \fill (0,0) node[anchor=north east]{$x$} circle (2pt);
    \fill (2,3.464) node[anchor=south] {$y$} circle (2pt);
    \fill (4,0) node[anchor=north west] {$z$} circle (2pt);
    \draw[->] (0,0) +(60:2pt) -- (60:4cm-2pt) ;
    \draw[->] (60:4cm) +(300:2pt) -- +(300:4cm-2pt);
\draw[<-,dashed] (2pt,0pt) --(4cm-2pt,0cm);
 \end{tikzpicture}
    \caption{Circular relation}
    \label{fig:circular}
\end{figure}

Garrett Birkhoff (1911--1996) gave an affirmative answer while studying the lattice of all equivalence relations~\cite[p.~446]{birkhoff}. The transitive property can be represented by a directed graph, as shown in Figure~\ref{fig:transitive}. The property says that if two edges are present, then a third edge (shown dashed) appears. 
Birkhoff\footnote{Garrett Birkhoff, who held the George Putnam Professorship of Mathematics at Harvard, should not be confused with his father, the mathematician George David Birkhoff (1884--1944), who was the fifth holder of the Perkins chair.} reversed the direction of the dashed edge, as shown in Figure~\ref{fig:circular}, to produce a property that he called ``circularity.'' This property of a relation~\(\sim\) says that if \(x\sim y\) and \(y\sim z\), then \(z\sim x\). Birkhoff presumably was unaware that 
Vailati had introduced circular relations under the name ``reversive'' some four decades earlier 
\cite[p.~162]{vailati}.

A reflexive and circular relation is automatically symmetric. Indeed, if \(x\sim y\), then specializing circularity to the case that \(z\) equals~\(y\) and invoking 
reflexivity to say that \(y\sim y\) yields that \(y\sim x\). In the presence of symmetry, the properties of transitivity and circularity are identical, so the following characterization is workable.

\begin{definition}[Birkhoff]
    An equivalence relation is a relation that is reflexive and circular.
    \label{def:Birkhoff}
\end{definition}

The value of this definition is not  solely the brevity obtained by reducing three hypotheses to two. More significant is the demonstration that symmetry is a subsidiary property,
not a crucial element of the \emph{definition} of equivalence relation. The essence of the concept of equivalence lies elsewhere. 

Birkhoff apparently felt no compulsion to convert colleagues to his definition. The groundbreaking \emph{Survey of Modern Algebra},\footnote{The first edition appeared in 1941, and the book is still available in a fifth edition \cite{birkhoffmaclane}.} which he coauthored with Saunders Mac Lane (1909--2005), states Hasse's definition of equivalence relation and relegates the notion of circularity to an exercise. 
I regret that I was too naive to ask Birkhoff about this circumstance when he was my designated faculty advisor during my first semester of college. Birkhoff and Mac Lane's later reminiscences \cite{reminiscences} celebrating the fiftieth anniversary of their book are silent about circularity, as is Mac Lane's autobiography \cite{maclane}. 

Definition~\ref{def:Birkhoff} does require reflexivity as a hypothesis. Indeed, the relation on the set \(\{x,y,z\}\) consisting of the three ordered pairs \( (x,y) \), \((y,z)\), and \((z,x)\) is circular, has domain equal to the underlying set, but is neither reflexive nor symmetric. 

\begin{figure}
    \centering
 \begin{tikzpicture}[>=Stealth,thick]
    \fill (0,0) node[anchor=north east]{$x$} circle (2pt);
    \fill (0,0) +(60:4cm) node[anchor=south] {$y$} circle (2pt);
    \fill (4,0) node[anchor=north west] {$z$} circle (2pt);
      \draw[->] (0,0) +(60:2pt) -- (60:4cm-2pt) ;
    \draw[->] (4,0) +(120:2pt) -- +(120:4cm-2pt);
    \draw[<-,dashed] (2pt,0pt) --(4cm-2pt,0cm); 
 \end{tikzpicture}
    \caption{Euclidean relation}
    \label{fig:Euclidean}
\end{figure}
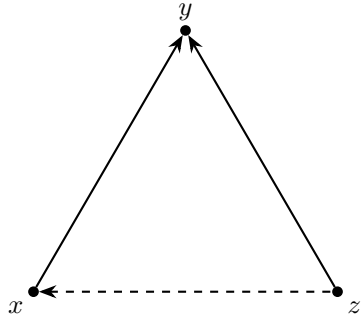

The two questionable aspects of Definition~\ref{def:hasse} can be addressed simultaneously by reversing a second arrow in the directed graph (as shown in Figure~\ref{fig:Euclidean}). The Euclidean property of a relation~\(\sim\) says that if \(x\sim y\) and \(z\sim y\), then \(z\sim x\). (Since the letters are implicitly universally quantified, interchanging \(x\) and~\(z\) shows that \(x\sim z\) too.) 

A reflexive and Euclidean relation is automatically symmetric. (Replace \(z\) by~\(y\) in Figure~\ref{fig:Euclidean}.) 
Therefore Definition~\ref{def:Birkhoff} admits a parallel formulation, which was known to De~Amicis \cite{deamicis} (whose adjective for a Euclidean relation was ``adequative''). 

\begin{definition}
    An equivalence relation is a relation that is reflexive and Euclidean.
\end{definition}

This statement is not the last word. 
Set \(z\) equal to~\(x\) in Figure~\ref{fig:Euclidean} to see that if \(\sim\)~is a Euclidean relation, and \(x\) belongs to the domain of the relation, then \(x\sim x\). In other words, a Euclidean relation is reflexive on its domain. This remark justifies the ultimate variation of the definition.

\begin{definition}
\label{def:optimal}
An equivalence relation on a set 
is a Euclidean relation whose domain is the whole set. 
\end{definition}

Although this final definition is implicit in the 1892 paper of De Amicis \cite[note~5]{deamicis}, the first explicit instance I know in the refereed literature is a 1967 note by Charles Buck \cite{buck}, who was active in mathematics education during the ``New Math'' era. He stated that the definition appeared in pedagogical materials for high-school students in 1965. Taking  Definition~\ref{def:optimal} as a starting point, one obtains all the standard properties of equivalence relations as a theorem.

\begin{theorem}
An equivalence relation is reflexive, symmetric, transitive, and circular. 
\end{theorem}

Thus Euclid's first common notion can be regarded as the fundamental property characterizing equivalence relations. 
This heuristic seems to have been folklore around 1870, before the concept of equivalence relation had a formal definition. 
Both the celebrated Dedekind and the unsung Charles M\'eray (1835--1911) chose to verify the Euclidean property of the relevant relation when justifying the existence of what are now called equivalence classes \cite{dirichlet, meray}. 

Definition~\ref{def:optimal} has minimal hypotheses, but is the statement pedagogically effective? 
I think so, because students routinely encounter relations that are transitive without being equivalence relations. The inequality relation~\(\le\) clearly is transitive but not Euclidean; and the same is true for the subset relation~\(\subseteq\),  for the relation \(P \Rightarrow Q\) in propositional logic, and for the relation ``\(k\) divides~\(n\)'' in number theory. Such examples validate the assertion that the Euclidean property---not transitivity---is the decisive condition for testing whether a relation is an equivalence. 

Accordingly, the question arises of who first formulated the Euclidean property for arbitrary binary relations. This significant event cannot predate 1850, for no general theory of relations existed before De~Morgan's work. 

My candidate for the perspicacious pioneer is Marian Evans (1819--1880), who was a professional writer and an amateur mathematician. In a discussion in the early 1850s with the prominent polymath Herbert Spencer (1820--1903), she formulated the axiom, ``Things that have a constant relation to the same thing have a constant relation to each other'' 
\cite[p.~162]{spencer}. 
Traces of the constant interest that Evans had in mathematics can be found in the acclaimed novels she began writing a few years later under the pen name George Eliot \cite{ball}. 

\begin{figure}
    \centering
\includegraphics[width=4cm]{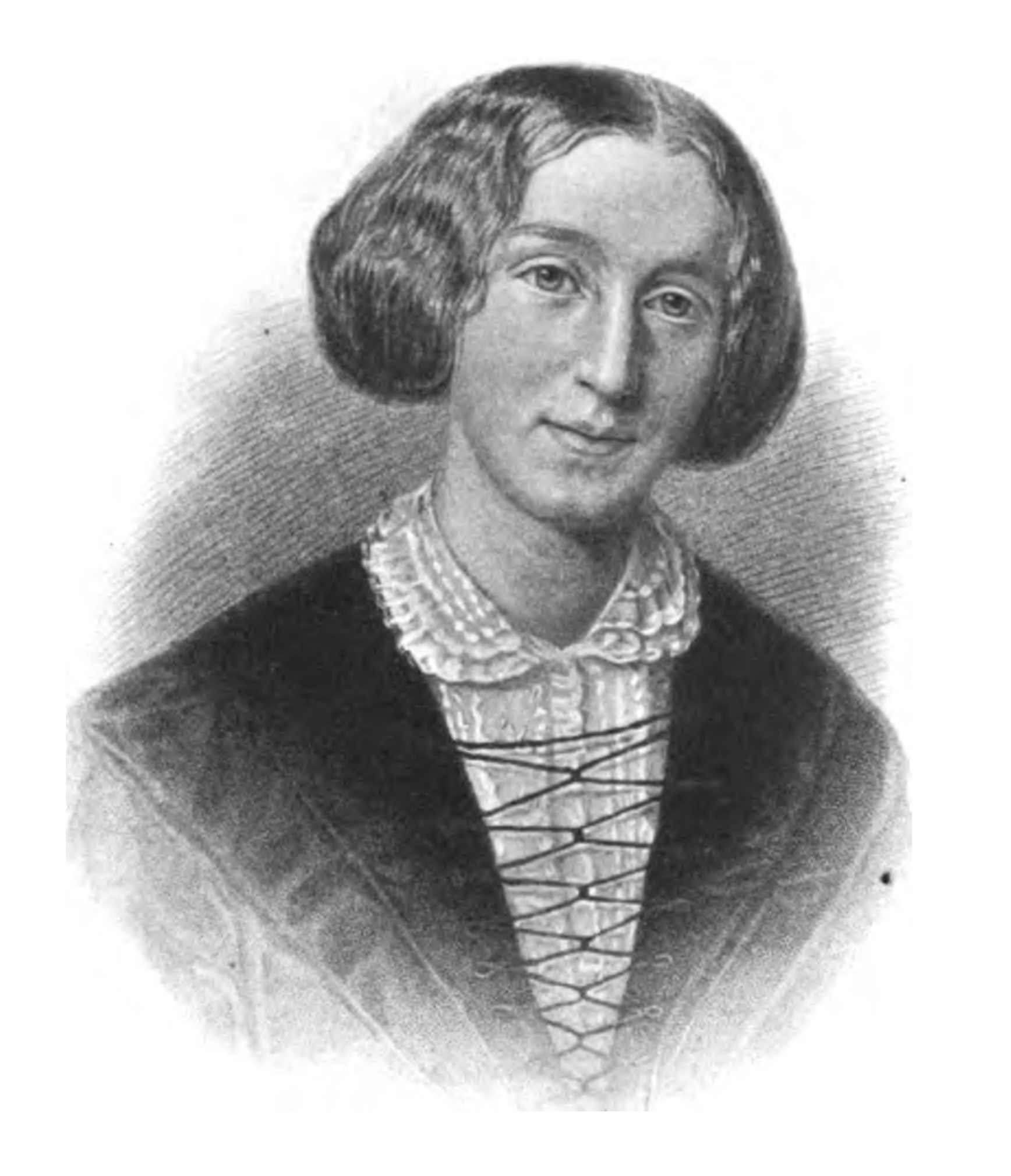}
    \caption{George Eliot in 1850 (public domain image \cite[frontispiece]{eliot})}
\end{figure}

\section*{One thing after another}
An early triumph of the nascent concept of equivalence relation was the publication in the 1870s of two constructions of the set of real numbers, starting from the set of rational numbers. 
Dedekind's method \cite{cuts} was to partition the rational numbers into two ordered segments (a partition of a set being a concrete way to realize an equivalence relation).  
Cantor's method, anticipated by M\'eray \cite{meray}, was to create the real numbers as equivalence classes of Cauchy sequences of rational numbers \cite{cantor-cauchy}. The latter procedure requires a definition of Cauchy sequence, which is the next topic.

The definition in formula~\eqref{eq:Cauchy} of the introduction displays an intimidating thicket of quantifiers. 
Even students who successfully parse the definition 
of ``\(\lim_{n\to\infty} x_n = L\),'' 
namely,
\begin{equation}
(\forall \varepsilon>0) (\exists N) (\forall m\ge N) ( |x_m-L|<\varepsilon),
\label{eq:limit}
\end{equation}
may balk at the fourth quantifier in~\eqref{eq:Cauchy}.

Yet the underlying idea is simple: a sequence can be viewed as an iterative process for approximating some number. A young child, unacquainted with trigonometry, has no difficulty finding a number that  equals its own cosine. The tot need only open a calculator app on a parent's phone and tap the \texttt{cos} button repeatedly until the display stops changing. The official definition of limit merely provides a guarantee that further taps will not change the result, within the accuracy of the calculator. 

Knowing the value of the limit~\(L\) in advance is not critical. If the goal is phrased as certifying that a suitable stopping point has been reached in an approximation scheme, then Cauchy sequences can be described by a statement parallel to~\eqref{eq:limit}, using only three quantifiers. 

\begin{definition}
\label{def:Heine}
A sequence \(\{x_n\}_{n=1}^\infty\) of real numbers is a Cauchy sequence when 
\begin{equation}
(\forall \varepsilon>0) (\exists N) (\forall m\ge N) ( |x_m-x_N|<\varepsilon).
\label{eq:Heine}
\end{equation}
\end{definition}

If this three-quantifier statement holds, then the triangle inequality implies that
\begin{equation*}
    |x_m-x_n| \le |x_m-x_N| + |x_n-x_N| < 2\varepsilon
    \qquad \text{when \(m\ge N\) and \(n\ge N\).}
\end{equation*}
Consequently, the four-quantifier statement~\eqref{eq:Cauchy} holds with the innocuous replacement of~\(\varepsilon\) by~\(2\varepsilon\) (allowable since \(\varepsilon\)~is universally quantified). 

Stated informally in words instead of symbols, the definition with four quantifiers says that eventually, the tail of the sequence has diameter as small as desired. The definition with three quantifiers says that eventually, the tail of the sequence lies inside a specific ball of radius as small as desired. The second statement is more concrete, hence easier for novices to learn. 
Students at my university would bless an instructor who discloses that once they have mastered the definition of limit, they also know the definition of Cauchy sequence: simply replace~\(L\) with the approximate value~\(x_N\). 

For half a century, I uncritically wrote the standard four quantifiers, so discovering Definition~\ref{def:Heine} was a shocking revelation. My amazed reaction echoed the exclamation of Thomas Hobbes (1588--1679) on his chance encounter with the Pythagorean theorem: ``this is impossible!''\ (related by John Aubrey \cite[p.~332]{aubrey}). 
Embarrassed by the gap in my knowledge of rudimentary analysis, I compiled the following history to fathom how the needlessly elaborate statement~\eqref{eq:Cauchy} became canonical. 

Mathematicians of the 18th century used infinite series to great effect, showing for example that
\begin{equation}
    \sum_{n=1}^\infty \frac{1}{n^2} = \frac{\pi^2}{6}, 
    \label{eq:Basel}
\end{equation}
yet the foundations of the subject were shaky. Indeed, when Jean le Rond d'Alembert (1717--1783) introduced a special case of the ratio test, he lamented that the proof 
of~\eqref{eq:Basel} given by
Johann Bernoulli (1667--1748) glossed over a ``mystery {\ldots} in the theory of series; a theory that appears to me still very imperfect'' \cite[p.~183]{alembert}. 

Unsurprisingly, the first steps toward a rigorous characterization of convergence were taken in the context of the sequence of partial sums of an infinite series. No formal theory of quantifiers was yet available, so the initial portrayals of convergence appear out of focus to modern eyes. Retouching the original work to fit current conventions requires some guesswork about the founders' viewpoints.

Jos\'e Anast\'acio da Cunha (1744--1787), who died two and a half years before Cauchy was born, defined an infinite series to be convergent when the sum can be terminated at some point after which adding an arbitrary number of additional terms makes a negligible change \cite[p.~106]{cunha}. 
This statement can plausibly be interpreted as a version of ``Cauchy's criterion'' for convergence \cite{queiro}, indeed the three-quantifier form~\eqref{eq:Heine}.

More explicit is the treatment 
in Cauchy's 1821 textbook, based on his lectures at the \'Ecole Polytechnique: A necessary and sufficient condition for convergence of an infinite series \(u_0 + u_1 + \cdots\) is that when \(n\)~increases, the sums 
\begin{equation}
\begin{aligned}
    & u_n +  u_{n+1} \\
    & u_n  + u_{n+1}+u_{n+2} \\
    & \qquad \vdots
\end{aligned}
\label{eq:tail}
\end{equation}
end up having absolute value less than any assigned (positive) bound \cite[pp.~124--126]{cauchy-cours}. I interpret Cauchy's statement as meaning that 
\begin{equation*}
    \lim_{n\to\infty} \sup_{p\ge 1} |u_n+\cdots + u_{n+p}|=0
\end{equation*}
(in modern notation), which corresponds to the four-quantifier condition~\eqref{eq:Cauchy}.

Bernard Bolzano (1781--1848) had 
already formulated Cauchy's criterion four years earlier,
stating as a theorem that a series converges when the indicated sums~\eqref{eq:tail} 
``stay smaller than any given quantity if one has taken \(n\) large enough'' \cite[\S~7]{bolzano}. The details of Bolzano's attempted (but failed) proof suggest to me that he meant \((\exists n)\) rather than \((\exists N)(\forall n\ge N)\). 

Thus the three-quantifier criterion~\eqref{eq:Heine} and the four-quantifier version~\eqref{eq:Cauchy} were not clearly distinguished when the concept of convergence was first formalized, just as one does not differentiate nowadays between \(\epsilon\) and~\(\varepsilon\) (the two forms of the letter epsilon). Of course, Bolzano\footnote{See also  accounts in this journal of Bolzano's work by Steve Russ \cite{russ} and by P.~Maritz \cite{maritz}.} and Cauchy and Cunha could not rigorously prove that a Cauchy sequence of real numbers converges, for a precise construction of the real numbers and a recognition of the need for a completeness axiom were decades in the future. 

The splendid abstract idea that a real number \emph{is} an equivalence class of Cauchy sequences of rational numbers goes back to M\'eray \cite[\S~IV]{meray}, who used the traditional language of Cauchy without explicit quantifiers. Cantor's independent treatment of this construction \cite{cantor-cauchy} explicitly supplied the \(\varepsilon\) and~\(N\) in condition~\eqref{eq:Cauchy}. Eduard Heine (1821--1881), Cantor's elder colleague at the University of Halle, wrote a polished exposition \cite{heine} of Cantor's method and stated the three-quantifier condition~\eqref{eq:Heine}.

Thus the twin characterizations \eqref{eq:Cauchy} and~\eqref{eq:Heine} not only hatched together but were still sharing a nest when Cauchy sequences became recognized as a species of mathematical objects in their own right. The two coexisted for another half century, but ultimately the ponderous definition~\eqref{eq:Cauchy} knocked the lighter sibling~\eqref{eq:Heine} out of the nest to expire from neglect. 

To account for the waning of the three-quantifier definition, I propose the analogy of constructed auxiliary languages, which were fashionable in the early 20th century (the idealistic goal being to facilitate international understanding). Peano developed \emph{Latino sine flexione}; Carnap, Maurice Fr\'echet (1878--1973), and M\'eray all were devotees of Esperanto; and Couturat was an instigator of \emph{Ido}, a child of Esperanto.  
Yet in the 21st century, the constructed languages of interest to  mathematicians are the ones used to communicate with machines for such purposes as scientific computing, modeling, typesetting, generating conjectures, and validating proofs. Fashions change.

The zenith of the sleek three-quantifier definition was the first decade of the 20th century. Notable backers included English mathematicians such as E.~T. Whittaker (1873--1956) in the first edition\footnote{Subsequent editions had G.~N. Watson (1886--1965) as coauthor.} of \emph{A Course of Modern Analysis} \cite[p.~8]{whittaker}, 
E.~W. Hobson (1856--1933) in his treatise on real functions \cite[p.~26]{hobson}, Jourdain \cite[p.~205]{jourdain-cauchy}, and T. J. I'A. Bromwich (1875--1929) in his book on infinite series \cite[p.~8]{bromwich}; French mathematicians such as \'Edouard Goursat (1858--1936) in his \emph{Cours d'analyse} \cite[p.~369]{goursat} and Fr\'echet in his thesis \cite[p.~23]{frechet} supervised by Jacques Hadamard (1865--1963); and the Belgian mathematician
Charles-Jean de~la Vall\'ee Poussin (1866--1962) in his \emph{Cours d'analyse} \cite[p.~8]{vallee}.
As the decade ended, the darkness of war was on the horizon, and the G\"ottingen lectures on real functions \cite[p.~92]{caratheodory} by the Greek--German mathematician Constantin Carath\'eodory (1873--1950) in the summer of 1914 were the last ray of sunshine 
for property~\eqref{eq:Heine}
before the long nightmare began.

The global conflict affected mathematicians in many ways \cite{guns}. Russell lost his lectureship for promulgating pacifist views unwelcome to the government (an echo of what befell Bolzano in his country a century earlier). Jourdain, who was a disciple of Russell's mathematics but not his politics, fervently wished to contribute to the war effort.
The Stephen Hawking of his time---a brilliant mind locked in a body ravaged by a progressive neurological disease---Jourdain was rebuffed by a  military unable to see past his disability.  
On the continent, Couturat was possibly the first French fatality of World War~I, killed on the very day that war was declared; he died in a car crash involving a military vehicle carrying urgent mobilization orders.  Fr\'echet applied his language skills to serve as an interpreter for the British army. Hadamard's two eldest sons died fighting against Germany in the endless Battle of Verdun.\footnote{See Jean-Pierre Kahane's article in this journal \cite{kahane} for more about Hadamard's life and legacy.} 
On the other side, Carath\'eodory's G\"ottingen colleague Richard Courant (1888--1972) was gravely wounded by British forces in the Battle of Loos. 

Courant had an outsized impact on the mathematical community in the postwar world \cite{reid}. Everyone who has learned from a Springer-Verlag book with a yellow cover has Courant to thank for persuading Ferdinand Springer (1881--1965) to jump into mathematics publishing.
Realizing a dream of Felix Klein (1849--1925), Courant built a mathematical institute in G\"ottingen. When the Nazi regime later nullified Courant's achievement by forcing out most of the faculty in 1933, Courant relocated to New York and created another renowned institute. His treatise about methods of mathematical physics, familiarly known as ``Courant--Hilbert,'' is still widely cited nearly a century after the initial publication.\footnote{For reminiscences of Courant published in this journal, see articles by Christopher~R. Friedrichs \cite{friedrichs} and David~E. Rowe \cite{rowe} as well as a transcription of Courant's own words \cite{courant-talk}.} 

The relevance of this history is that the tipping point for the definition of Cauchy sequence seems to be the 1927 first volume of Courant's calculus textbook, which happened to use the four-quantifier definition \cite[p.~29]{courant}. This influential volume went through multiple editions and reprintings both in German and in English translation.
Although the three-quantifier definition was still extant, appearing in the 1927 second edition of the book on set theory \cite[p.~103]{hausdorff2} by Felix Hausdorff (1868--1942) and in the 1933 first volume on topology \cite[p.~196]{kuratowski} by Kazimierz Kuratowski (1896--1980), its days were numbered.  
A quarter century later, a new generation of outstanding expositors adopted the four-quantifier definition in their books, such as 
\emph{Principles of Mathematical Analysis} \cite[p.~39]{rudin}
by Walter Rudin (1921--2010)
and \emph{Mathematical Analysis} \cite[p.~66]{apostol}
by Tom~M. Apostol (1923--2016).
 
Students of these gurus naturally copied the four-quantifier definition and passed it on to their own pupils. The books from which I learned, and the books from which I later taught, all used the four-quantifier definition. Will the simpler three-quantifier definition of Cauchy sequence rise again one day, like a phoenix? 

\subsection*{A remark on terminology}
Although I have been using the words ``Cauchy sequence'' throughout the discussion, this designation did not exist before the 20th century. No name was needed until the work of Cantor, whose own terminology was ``fundamental sequence'' \cite[p.~567]{cantor-fundamental}.
As far as I know, the 
appellation ``Cauchy sequence'' 
(\emph{suite de Cauchy}) first 
appeared in Fr\'echet's thesis \cite[p.~24]{frechet}. 

I have restricted for simplicity to the motivating case of a Cauchy sequence of real numbers, yet the concept makes sense whenever there is a notion of distance between points. The natural setting for Cauchy sequences is a metric space, which is the next subject.

\section*{So far away}
The concepts of an equivalence relation and a Cauchy sequence developed on a time scale of decades, as a trend toward abstraction took hold in 19th-century mathematics.
On the other hand, the notion of a metric space has a clear-cut 20th-century origin in Fr\'echet's thesis, which provided axioms equivalent to those stated in my introduction \cite[p.~30]{frechet}. Hausdorff\footnote{See articles in this journal by Walter Purkert \cite{purkert} and Charlotte K. Simmons \cite{simmons} for accounts of Hausdorff's life and work.} standardized axioms \eqref{eq:id}--\eqref{eq:tri} in his foundational book on set theory and introduced the terminology ``metric space'' \cite[p.~211]{hausdorff}, nomenclature that became widely accepted despite Fr\'echet's reservations \cite[p.~62]{frechet-book}. 

These standard axioms are not minimal! One minor simplification---sufficiently well known to be recorded in Wikipedia---is to drop the assumption that the distance function is nonnegative.
As long as the distance is real-valued, so that inequalities make sense, setting \(z\) equal to~\(x\) in the triangle inequality and invoking the coincidence axiom shows that
\begin{equation*}
0=d(x,x) \le d(x,y) + d(y,x),
\end{equation*}
so by the symmetry axiom, the quantity \(d(x,y)\) cannot be negative. The redundancy of the hypothesis of nonnegativity was pointed out in 1926 by Adolf Lindenbaum (1904--1941) while still a student \cite[p.~211]{lindenbaum}.

Moreover, I find the coincidence axiom~\eqref{eq:id} more intelligible when expressed in geometric language.
The property says that the zero set of the distance function~\(d\) is the diagonal, the set of ordered pairs having equal coordinates.

The characteristic feature of metric spaces is the triangle inequality
(already invoked in the previous section). The name harks back to Proposition~20 of Book~I of Euclid's \emph{Elements}: ``In any triangle two sides taken together in any manner are greater than the remaining one'' (Heath's translation \cite[p.~286]{heath}). Apparently, 
Fr\'echet wrote the triangle inequality in the form~\eqref{eq:tri} to ensure transitivity of the informally stated relation ``the distance between \(x\) and~\(y\) is infinitesimal.'' Since circularity is a more powerful property than transitivity, the triangle inequality can be strengthened to a circular version,
\begin{equation}
    d(z,x) \le d(x,y) + d(y,z).
    \label{eq:cir}
\end{equation}

Garrett Birkhoff \cite[p.~466]{birkhoff-metric} gave the circular inequality~\eqref{eq:cir} a pleasing real-world interpretation, which I paraphrase as follows. If you take a plane trip from city~\(x\) to city~\(z\) with a stop in city~\(y\), then your return trip will be shorter if you take a direct flight home from~\(z\) to~\(x\).

Combined with the coincidence axiom, the circular inequality implies that the distance function is symmetric. Indeed, setting \(z\) equal to~\(y\) in~\eqref{eq:cir} shows that \(d(y,x) \le d(x,y)\). The variables are  universally quantified, so the reverse inequality follows by switching \(x\) and~\(y\), whence \(d(x,y)=d(y,x)\). 
In the presence of symmetry, the two triangle inequalities \eqref{eq:tri} and~\eqref{eq:cir} agree.
Hence the number of axioms for a metric space can be reduced from three to two by replacing the triangle inequality~\eqref{eq:tri} with the circular version~\eqref{eq:cir} and omitting the symmetry axiom~\eqref{eq:sym}.

The Dutch mathematician 
David van Dantzig (1900--1959) anticipated Birkhoff by a decade in stating these two axioms for metric spaces \cite[p.~598]{dantzig}, at the same time citing Lindenbaum's prior publication, which contains the circular inequality in a footnote \cite[p.~211]{lindenbaum}.  Lindenbaum wrote that he got the circular inequality from Piotr Szyma\'nski (1900--1965), a fellow student at the University of Warsaw. 

Lindenbaum's own innovation was to replace inequality~\eqref{eq:tri} not by~\eqref{eq:cir} but
by an inequality corresponding to Figure~\ref{fig:Euclidean}:
\begin{equation}
d(x,z) \le d(y,x) + d(y,z). 
\label{eq:ell}
\end{equation}
I will call~\eqref{eq:ell} the elliptical inequality, interpreting it as saying that the sum of the distances from a point~\(y\) on an ellipse to the foci \(x\) and~\(z\) exceeds the distance between the foci. 

Like the circular inequality, 
the elliptical inequality forces symmetry by setting \(z\) equal to~\(y\). Since the right-hand side of~\eqref{eq:ell} is unchanged when \(x\) and~\(z\) are switched, the left-hand side can be written either as \(d(x,z)\) or as \(d(z,x)\). A parallel argument shows that the right-hand side can equally well be \(d(x,y) + d(z,y)\). This flexibility 
means that if you take the triangle in Figure~\ref{fig:Euclidean} and reassign the directions of the three arrows at random, you will obtain an elliptical inequality half the time, and the other cases will be split between the circular inequality and the standard triangle inequality. Hence Lindenbaum's definition of a metric space seems the most natural. 

\begin{definition}[Lindenbaum]
A metric space is a set~\(E\) together with a real-valued function~\(d\) on \(E\times E\) such that 
\begin{enumerate}
    \item the zero set of~\(d\) is the diagonal, and 
    \item the distance function \(d\) satisfies the elliptical inequality~\eqref{eq:ell}. 
\end{enumerate}
\end{definition}

\begin{theorem}
In a metric space, the distance function is nonnegative and symmetric. Moreover, the triangle inequality~\eqref{eq:tri} and the circular inequality~\eqref{eq:cir} hold. 
\end{theorem}

Lindenbaum is remembered  mainly for contributions to logic, in particular his collaboration with Alfred Tarski\footnote{For an appreciation of Tarski's life and work, see an article in this journal by Steven~R. Givant \cite{givant}.} (1901--1983), but Polish mathematicians of the time were well aware of Lindenbaum's
clever definition of metric space. Stefan Banach (1892--1945) stated the Fr\'echet--Hausdorff axioms in his seminal book on functional analysis but mentioned Lindenbaum's improvement in a footnote \cite[p.~8]{banach}, as did Lindenbaum's thesis supervisor Wac{\l}aw Sierpi\'nski (1882--1969) in his book on point-set topology \cite[p.~75]{sierpinski}. Kuratowski (who also had taught Lindenbaum) jettisoned the traditional axioms in favor of Lindenbaum's \cite[pp.~82--83]{kuratowski}.

In America, Solomon Lefschetz (1884--1972) adopted Lindenbaum's definition of metric space in a book \cite[p.~5]{lefschetz} that standardized the word ``topology'' for the subject previously known as ``analysis situs.'' The indirect participation of Alice Berg Hayes (1888--1978) should be acknowledged. She and Lefschetz met as graduate students, and she had a teaching job before they married. Although the social mores of the era precluded her from contemplating an academic career, she made a vital contribution to mathematics through helping her husband cope with the physical and psychological challenges of his disability. (In his short-lived first career as an engineer, Lefschetz had lost both hands in the explosion of a high-voltage transformer that he was testing.) 

Although Lindenbaum's definition was noticed by some prominent researchers, most authors of introductory textbooks failed to follow suit, one
exception being E.~T. Copson (1901--1980) in a Cambridge Tract \cite[p.~21]{copson}. I~explain this circumstance through a metaphor.

Fr\'echet and Hausdorff constructed the initial level of an avant-garde skyscraper, the headquarters for two flourishing new corporations. The first was point-set topology in general abstract spaces, which had \emph{Fundamenta Mathematicae} as house journal. Sierpi\'nski was one of the founding editors in 1920, and Kuratowski joined the board in 1929.  
The other new entity was functional analysis, whose name was provided in 
a book \cite{levy} by Paul L\'evy (1886--1971), taking a hint from his and Fr\'echet's teacher Hadamard.
The supporting periodical was \emph{Studia Mathematica}, of which Banach was a founding editor. 
Two decades after the groundbreaking, when Lindenbaum proposed an elegant upgrade to the front entrance of the edifice, crews were too busy building the upper stories to pay much attention. Perhaps Lindenbaum's new design would eventually have been exhumed from the archives if not for the man-made disaster that shut down the construction site for half a dozen years.

This catastrophe was the second world war, antithetical to rational thought. Amid the
universal suffering, some  characters in this chronicle found abstract mathematics to be a solace in the dark, an affirmation of faith in the human spirit, or, in extremis, an agency for preserving if not one's life then one's soul. 

When war broke out, Tarski fortuitously was overseas at Quine's invitation.
Returning home to Warsaw was out of the question, but
Tarski was aided by Carnap, Quine, and Russell to obtain an immigrant visa. The planned International Congress of Mathematicians, scheduled for September 1940 at Harvard, had to be canceled. The German forces occupying Poland ransacked the universities, but Kuratowski got help from the Swiss consulate to preserve the draft second volume of his topology book by sending it to Geneva, ``as one throws things of value, if only to oneself, at the last moment out of a burning house.''\footnote{The time-shifted quotation is from Lord Dunsany (1878--1957) in the previous war \cite[preface]{dunsany}.} 
Rudin's family had got out of Austria after the Anschluss, going first to Switzerland and then to France. After the collapse of the French Third Republic, Rudin fled to England.
In Vichy France, the peril for Jews forced L\'evy literally to head for the hills and go into hiding. 
Hadamard was able to secure a visiting position in America, but his last son died fighting with the Free French in North Africa. Tragically, many could find no refuge. Hausdorff and his wife, Tarski's parents, and 
Lindenbaum and his wife, who was the logician Janina Hosiasson (1899--1942), all perished in the Holocaust. 

\section*{Conclusion}
Mathematical vocabulary evolves over time. Cauchy could call a list of numbers either a ``sequence'' or a ``series,'' but today's calculus teachers tell students that confusing these words is like mixing up American football with the football played in the rest of the world. Russell could speak of ``definition by means of the principle of abstraction,'' but this terminology has been supplanted by the concepts of equivalence classes and quotient structures. The generalizations of metric spaces called ``class~\((V)\)'' in Fr\'echet's thesis are obsolete, for they turned out to be metrizable \cite{chittenden}. 

My leitmotif is that mathematical  language can reach equilibrium by happenstance instead of by design. The steady state need not be a global extremum. 

This article has discussed the history of three definitions whose standard forms illuminate their subject with moonlight rather than sunlight. 
The takeaway message is that alternative definitions shine  brighter and sharper.   
Equivalence relations can be characterized by \emph{one} condition, the Euclidean property. Metric spaces can be defined by \emph{two} axioms: the elliptical triangle inequality holds, and the zero set of the distance function equals the diagonal. Cauchy sequences can be specified through an inequality governed by \emph{three} quantifiers. It's as simple as one, two, three.

\section*{Declaration} 
No funding was received to assist with the preparation of this manuscript.

\bibliographystyle{abbrvurl}
\bibliography{trois}

\begin{thebibliography}{10}

\bibitem{apostol}
T.~M. Apostol.
\newblock {\em Mathematical Analysis: {A} Modern Approach to Advanced Calculus}.
\newblock Addison-Wesley, Reading, MA, first edition, 1957.

\bibitem{asghari}
A.~Asghari.
\newblock Equivalence: an attempt at a history of the idea.
\newblock {\em Synthese}, 196(11):4657--4677, 2019.
\newblock \href {https://doi.org/10.1007/s11229-018-1674-2} {\path{doi:10.1007/s11229-018-1674-2}}.

\bibitem{guns}
D.~Aubin and C.~Goldstein, editors.
\newblock {\em The War of Guns and Mathematics}.
\newblock Amer. Math. Soc., Providence, 2014.
\newblock \href {https://doi.org/10.1090/hmath/042} {\path{doi:10.1090/hmath/042}}.

\bibitem{ball}
D.~Ball.
\newblock Mathematical contrariness in {G}eorge {E}liot's novels.
\newblock In R.~Tubbs, A.~Jenkins, and N.~Engelhardt, editors, {\em Palgrave Handbook of Literature and Mathematics}, chapter~6, pages 97--111. Palgrave Macmillan, Cham, 2021.
\newblock \href {https://doi.org/10.1007/978-3-030-55478-1_6} {\path{doi:10.1007/978-3-030-55478-1_6}}.

\bibitem{banach}
S.~Banach.
\newblock {\em Th\'eorie des {o}p\'erations {l}in\'eaires}.
\newblock Garasi\'nski, Warsaw, first edition, 1932.
\newblock URL: \url{http://kielich.amu.edu.pl/Stefan_Banach/e-operations.html}.

\bibitem{birkhoff}
G.~Birkhoff.
\newblock On the structure of abstract algebras.
\newblock {\em Proc. Camb. Philos. Soc.}, 31(4):433--454, 1935.
\newblock \href {https://doi.org/10.1017/S0305004100013463} {\path{doi:10.1017/S0305004100013463}}.

\bibitem{birkhoff-metric}
G.~Birkhoff.
\newblock Metric foundations of geometry. {I}.
\newblock {\em Trans. Amer. Math. Soc.}, 55:465--492, 1944.
\newblock \href {https://doi.org/10.2307/1990304} {\path{doi:10.2307/1990304}}.

\bibitem{reminiscences}
G.~Birkhoff and S.~Mac~Lane.
\newblock \emph{{A} {S}urvey of {M}odern {A}lgebra}: {T}he fiftieth anniversary of its publication.
\newblock {\em Math. Intelligencer}, 14(1):26--31, 1992.
\newblock \href {https://doi.org/10.1007/BF03024138} {\path{doi:10.1007/BF03024138}}.

\bibitem{birkhoffmaclane}
G.~Birkhoff and S.~Mac~Lane.
\newblock {\em A Survey of Modern Algebra}.
\newblock CRC Press, Boca Raton, fifth edition, 2010.
\newblock \href {https://doi.org/10.1201/9781315275499} {\path{doi:10.1201/9781315275499}}.

\bibitem{bisht}
R.~K. Bisht and R.~Wilson.
\newblock Bertrand {R}ussell.
\newblock {\em Math. Intelligencer}, 45(4):392--393, 2023.
\newblock \href {https://doi.org/10.1007/s00283-023-10300-7} {\path{doi:10.1007/s00283-023-10300-7}}.

\bibitem{bolzano}
B.~Bolzano.
\newblock {\em Rein {a}nalytischer Beweis des Lehrsatzes {da\ss} {z}wischen {je} {zwey} Werthen, die {ein} {e}ntgegengesetzetes Resultat {g}ew\"ahren, {w}enigstens {eine} {r}eelle Wurzel der Gleichung {l}iege}.
\newblock Haase, Prague, 1817.
\newblock URL: \url{https://hdl.handle.net/10338.dmlcz/400019}.

\bibitem{bourbaki}
N.~Bourbaki.
\newblock {\em Th\'eorie des {e}nsembles}.
\newblock Springer, Berlin, 2006.
\newblock Reprint of the 1970 edition.
\newblock \href {https://doi.org/10.1007/978-3-540-34035-5} {\path{doi:10.1007/978-3-540-34035-5}}.

\bibitem{bromwich}
T.~J.~I. Bromwich.
\newblock {\em An Introduction to the Theory of Infinite Series}.
\newblock Macmillan, London, first edition, 1908.
\newblock URL: \url{https://books.google.com/books?id=ZY45AAAAMAAJ}.

\bibitem{buck}
C.~Buck.
\newblock An alternative definition for equivalence relations.
\newblock {\em Math. Teacher}, 60:124--125, 1967.
\newblock URL: \url{https://www.jstor.org/stable/27957510}.

\bibitem{cantor-cauchy}
G.~Cantor.
\newblock Ueber die {A}usdehnung eines {S}atzes aus der {T}heorie der trigonometrische {R}eihen.
\newblock {\em Math. Ann.}, 5:123--132, 1872.
\newblock \href {https://doi.org/10.1007/BF01446327} {\path{doi:10.1007/BF01446327}}.

\bibitem{cantor}
G.~Cantor.
\newblock Ein {B}eitrag zur {M}annigfaltigkeitslehre.
\newblock {\em J. reine angew. Math.}, 84:242--259, 1877.
\newblock URL: \url{https://eudml.org/doc/148353}.

\bibitem{cantor-fundamental}
G.~Cantor.
\newblock Ueber unendliche, lineare {P}unktmannichfaltigkeiten.
\newblock {\em Math. Ann.}, 21:545--591, 1883.
\newblock \href {https://doi.org/10.1007/BF01446819} {\path{doi:10.1007/BF01446819}}.

\bibitem{caratheodory}
C.~Carath{\'e}odory.
\newblock {\em Vorlesungen {{\"u}}ber {r}eelle {Funktionen}}.
\newblock Teubner, Leipzig, first edition, 1918.
\newblock URL: \url{https://books.google.com/books?id=T4nMBa1fWo0C}.

\bibitem{carnap}
R.~Carnap.
\newblock {\em Abriss der Logistik}.
\newblock Springer, Vienna, 1929.
\newblock URL: \url{https://books.google.com/books?id=BloaAAAAIAAJ}.

\bibitem{cauchy-quotient}
A.~Cauchy.
\newblock M\'emoire sur une nouvelle th\'eorie des imaginaires, et sur les racines symboliques des \'equations et des \'equivalences.
\newblock {\em C. R. hebd. s{\'e}ances Acad. sci.}, 24:1120--1130, 1847.
\newblock URL: \url{https://gallica.bnf.fr/ark:/12148/bpt6k29812}.

\bibitem{cauchy-cours}
A.-L. Cauchy.
\newblock {\em Cours {d}'analyse de l'\'Ecole {p}olytechnique: Analyse {a}lg\'ebrique}.
\newblock Debure, Paris, 1821.
\newblock URL: \url{https://gallica.bnf.fr/ark:/12148/btv1b8626657t}.

\bibitem{chittenden}
E.~W. Chittenden.
\newblock On the equivalence of {\'e}cart and voisinage.
\newblock {\em Trans. Amer. Math. Soc.}, 18:161--166, 1917.
\newblock \href {https://doi.org/10.2307/1988857} {\path{doi:10.2307/1988857}}.

\bibitem{aubrey}
A.~Clark, editor.
\newblock {\em Aubrey's `Brief Lives'}, volume~I.
\newblock Oxford Univ. Press, London, 1898.
\newblock URL: \url{https://www.google.com/books/edition/_/WYtCAAAAIAAJ}.

\bibitem{copson}
E.~T. Copson.
\newblock {\em Metric Spaces}.
\newblock Univ. Press, Cambridge, 1968.
\newblock \href {https://doi.org/10.1017/CBO9780511566141} {\path{doi:10.1017/CBO9780511566141}}.

\bibitem{courant}
R.~Courant.
\newblock {\em Vorlesungen \"uber Differential- und Integralrechnung}, volume~1.
\newblock Springer, Berlin, first edition, 1927.
\newblock URL: \url{https://hdl.handle.net/2027/wu.89062908009}.

\bibitem{courant-talk}
R.~Courant.
\newblock Reminiscences from {H}ilbert's {G}\"ottingen.
\newblock {\em Math. Intelligencer}, 3(4):154--164, 1980/81.
\newblock \href {https://doi.org/10.1007/BF03022974} {\path{doi:10.1007/BF03022974}}.

\bibitem{couturat}
L.~Couturat.
\newblock {\em Les {p}rincipes des {m}ath\'ematiques}.
\newblock Alcan, Paris, 1905.
\newblock URL: \url{https://books.google.com/books?id=h7rqjsm_aw0C}.

\bibitem{cunha}
J.~A. {Da Cunha}.
\newblock {\em Principios {m}athematicos}.
\newblock Galhardo, Lisbon, 1790.
\newblock URL: \url{https://resolver.sub.uni-goettingen.de/purl?PPN590888331}.

\bibitem{alembert}
D'Alembert.
\newblock {\em Opuscules {m}ath\'ematiques}, volume~V.
\newblock Briasson, Paris, 1768.
\newblock URL: \url{https://gallica.bnf.fr/ark:/12148/bpt6k62424s}.

\bibitem{deamicis}
E.~De~Amicis.
\newblock Dipendenza fra alcune propriet\`a notevoli delle relazioni fra enti di un medesimo sistema.
\newblock {\em Riv. Mat.}, II:113--127, 1892.
\newblock URL: \url{https://books.google.com/books?id=PQ4MAAAAYAAJ}.

\bibitem{vallee}
C.-J. {de la Vall{\'e}e Poussin}.
\newblock {\em Cours {d}'analyse {i}nfinit{\'e}simale}, volume~I.
\newblock Gauthier-Villars, Paris, first edition, 1903.
\newblock URL: \url{https://books.google.com/books?id=r8rNAAAAMAAJ}.

\bibitem{DeMorgan}
A.~De~Morgan.
\newblock On the symbols of logic, the theory of the syllogism, and in particular of the copula, and the application of the theory of probabilities to some questions of evidence.
\newblock {\em Trans. Camb. Philos. Soc.}, 9 part I:79--127, 1856.
\newblock URL: \url{https://hdl.handle.net/2027/mdp.39015012112531}.

\bibitem{cuts}
R.~Dedekind.
\newblock {\em Stetigkeit und {i}rrationale {Zahlen}}.
\newblock Vieweg, Braunschweig, 1872.
\newblock URL: \url{https://books.google.com/books?id=n-43AAAAMAAJ}.

\bibitem{carroll}
C.~L. Dodgson.
\newblock {\em Euclid and His Modern Rivals}.
\newblock Macmillan, London, second edition, 1885.
\newblock URL: \url{https://hdl.handle.net/2027/coo.31924060288804}.

\bibitem{dunsany}
Dunsany.
\newblock {\em The Last Book of Wonder}.
\newblock Luce, Boston, 1916.
\newblock URL: \url{https://www.loc.gov/item/17002700/}.

\bibitem{eliot}
G.~Eliot.
\newblock {\em The Works of George Eliot}, volume XII: Poems.
\newblock Wheeler Publishing, New York, 1900.
\newblock URL: \url{https://www.google.com/books/edition/The_Works_of_George_Eliot_Poems/Ec8OAAAAIAAJ}.

\bibitem{frechet}
M.~Fr{\'e}chet.
\newblock Sur quelques points du calcul fonctionnel.
\newblock {\em Rend. Circ. Matem. Palermo}, 22:1--74, 1906.
\newblock \href {https://doi.org/10.1007/BF03018603} {\path{doi:10.1007/BF03018603}}.

\bibitem{frechet-book}
M.~Fr\'echet.
\newblock {\em Les {e}spaces {a}bstraits}.
\newblock Gauthier-Villars, Paris, 1928.
\newblock URL: \url{https://books.google.com/books?id=g_3uAAAAMAAJ}.

\bibitem{friedrichs}
C.~R. Friedrichs.
\newblock The {C}ourant circle as an extended family: {N}ew {R}ochelle and beyond.
\newblock {\em Math. Intelligencer}, 37(1):41--44, 2015.
\newblock \href {https://doi.org/10.1007/s00283-014-9523-8} {\path{doi:10.1007/s00283-014-9523-8}}.

\bibitem{gauss}
C.~F. Gauss.
\newblock {\em Disquisitiones {a}rithmeticae}.
\newblock Gerhard Fleischer, Leipzig, 1801.
\newblock URL: \url{https://resolver.sub.uni-goettingen.de/purl?PPN235993352}.

\bibitem{givant}
S.~R. Givant.
\newblock A portrait of {A}lfred {T}arski.
\newblock {\em Math. Intelligencer}, 13(3):16--32, 1991.
\newblock \href {https://doi.org/10.1007/BF03023831} {\path{doi:10.1007/BF03023831}}.

\bibitem{goursat}
{\'E}.~Goursat.
\newblock {\em Cours {d}'analyse {m}ath{\'e}matique}, volume~I.
\newblock Gauthier-Villars, Paris, first edition, 1902.
\newblock URL: \url{https://books.google.com/books?id=IRbvAAAAMAAJ}.

\bibitem{hartshorne}
R.~Hartshorne.
\newblock {\em Geometry: {E}uclid and Beyond}.
\newblock Springer, New York, 2000.
\newblock \href {https://doi.org/10.1007/978-0-387-22676-7} {\path{doi:10.1007/978-0-387-22676-7}}.

\bibitem{hasse}
H.~Hasse.
\newblock {\em H{\"o}here {Algebra}. {I}: {Lineare} {Gleichungen}}.
\newblock de Gruyter, Berlin, first edition, 1926.
\newblock URL: \url{https://hdl.handle.net/2027/mdp.39015035207730}.

\bibitem{hausdorff}
F.~Hausdorff.
\newblock {\em Grundz\"uge der Mengenlehre}.
\newblock von Veit, Leipzig, first edition, 1914.
\newblock URL: \url{https://books.google.com/books?id=KTs4AAAAMAAJ}.

\bibitem{hausdorff2}
F.~Hausdorff.
\newblock {\em Mengenlehre}.
\newblock de Gruyter, Berlin, second edition, 1927.
\newblock URL: \url{https://hdl.handle.net/2027/uc1.b4248809}.

\bibitem{heath}
T.~L. Heath, editor.
\newblock {\em The Thirteen Books of {E}uclid's {E}lements}, volume~1.
\newblock Univ. Press, Cambridge, first edition, 1908.
\newblock URL: \url{https://books.google.com/books?id=dkk6AQAAMAAJ}.

\bibitem{heine}
E.~Heine.
\newblock Die {E}lemente der {F}unctionenlehre.
\newblock {\em J. reine angew. Math.}, 74:172--188, 1872.
\newblock \href {https://doi.org/10.1515/crll.1872.74.172} {\path{doi:10.1515/crll.1872.74.172}}.

\bibitem{hobson}
E.~W. Hobson.
\newblock {\em The Theory of Functions of a Real Variable and the Theory of {Fourier}'s Series}.
\newblock Univ. Press, Cambridge, first edition, 1907.
\newblock URL: \url{https://name.umdl.umich.edu/ACM2112.0001.001}.

\bibitem{jourdain-cauchy}
P.~E.~B. Jourdain.
\newblock The introduction of irrational numbers.
\newblock {\em Math. Gaz.}, 4:201--209, 1908.
\newblock \href {https://doi.org/10.2307/3602961} {\path{doi:10.2307/3602961}}.

\bibitem{jourdain}
P.~E.~B. Jourdain.
\newblock On isoid relations and theories of irrational number.
\newblock In E.~W. Hobson and A.~E.~H. Love, editors, {\em Proceedings of the Fifth International Congress of Mathematicians}, volume~II, pages 492--496, Cambridge, 1913. Univ. Press.
\newblock URL: \url{https://www.mathunion.org/fileadmin/ICM/Proceedings/ICM1912.2/ICM1912.2.ocr.pdf}.

\bibitem{kahane}
J.-P. Kahane.
\newblock Jacques {H}adamard.
\newblock {\em Math. Intelligencer}, 13(1):23--29, 1991.
\newblock \href {https://doi.org/10.1007/BF03024068} {\path{doi:10.1007/BF03024068}}.

\bibitem{kuratowski}
C.~Kuratowski.
\newblock {\em Topologie I}.
\newblock Garasi\'nski, Warsaw, first edition, 1933.
\newblock URL: \url{http://pldml.icm.edu.pl/pldml/element/bwmeta1.element.zamlynska-9ee68a27-f16a-4074-a332-1d8b3b83f2a7}.

\bibitem{lefschetz}
S.~Lefschetz.
\newblock {\em Topology}.
\newblock Amer. Math. Soc., New York, 1930.
\newblock URL: \url{https://bookstore.ams.org/coll-12}.

\bibitem{legendre}
A.-M. Legendre.
\newblock {\em \'El\'ements de {g}\'eom\'etrie}.
\newblock Firmin Didot, Paris, first edition, 1794.
\newblock URL: \url{https://gallica.bnf.fr/ark:/12148/bpt6k1521831j}.

\bibitem{dirichlet}
P.~G. {Lejeune Dirichlet}.
\newblock {\em Vorlesungen \"uber {Z}ahlentheorie}.
\newblock Vieweg, Braunschweig, second edition, 1871.
\newblock Edited with additions by R.~Dedekind.
\newblock URL: \url{https://resolver.sub.uni-goettingen.de/purl?PPN30976923X}.

\bibitem{levy}
P.~L\'evy.
\newblock {\em Le\c{c}ons {d}'analyse {f}onctionnelle}.
\newblock Gauthier-Villars, Paris, 1922.
\newblock URL: \url{https://books.google.com/books?id=7TAPAAAAIAAJ}.

\bibitem{lindenbaum}
A.~Lindenbaum.
\newblock Contributions {\`a} l'{\'e}tude de l'espace m{\'e}trique {I}.
\newblock {\em Fund. Math.}, 8:209--222, 1926.
\newblock URL: \url{https://eudml.org/doc/214867}.

\bibitem{lucic}
Z.~Lu\v{c}i\'c.
\newblock Who proved {P}ythagoras's theorem?
\newblock {\em Math. Intelligencer}, 44(4):373--381, 2022.
\newblock \href {https://doi.org/10.1007/s00283-022-10205-x} {\path{doi:10.1007/s00283-022-10205-x}}.

\bibitem{maclane}
S.~Mac~Lane.
\newblock {\em Saunders {M}ac {L}ane---A Mathematical Autobiography}.
\newblock A K Peters, New York, 2005.
\newblock With a preface by David Eisenbud.
\newblock \href {https://doi.org/10.1201/9781439863640} {\path{doi:10.1201/9781439863640}}.

\bibitem{maritz}
P.~Maritz.
\newblock The {B}olzano house in {P}rague.
\newblock {\em Math. Intelligencer}, 23(2):52--55, 2001.
\newblock \href {https://doi.org/10.1007/BF03026628} {\path{doi:10.1007/BF03026628}}.

\bibitem{meray}
C.~M\'eray.
\newblock Remarques sur la nature des quantit\'es d\'efinies par la condition de servir de limites \`a des variables donn\'ees.
\newblock {\em Rev. Soc. sav. Sci. math. phys. nat. (2)}, 4:280--289, 1869.
\newblock URL: \url{https://gallica.bnf.fr/ark:/12148/bpt6k2026062}.

\bibitem{palais}
B.~Palais.
\newblock \(\pi\) is wrong!
\newblock {\em Math. Intelligencer}, 23:7--8, 2001.
\newblock \href {https://doi.org/10.1007/BF03026846} {\path{doi:10.1007/BF03026846}}.

\bibitem{calcolo}
G.~Peano.
\newblock {\em Calcolo Geometrico}.
\newblock Bocca, Turin, 1888.
\newblock URL: \url{https://books.google.com/books?id=5LJi3dxLzuwC}.

\bibitem{geometria}
G.~Peano.
\newblock {\em I {p}rincipii {di} Geometria {l}ogicamente {e}sposti}.
\newblock Bocca, Turin, 1889.
\newblock URL: \url{https://books.google.com/books?id=q14LAAAAYAAJ}.

\bibitem{arithmetices}
I.~Peano.
\newblock {\em Arithmetices {p}rincipia}.
\newblock Bocca, Turin, 1889.
\newblock URL: \url{https://books.google.com/books?id=UUFtAAAAMAAJ}.

\bibitem{peirce}
C.~S. Peirce.
\newblock Description of a notation for the logic of relatives, resulting from an amplification of the conceptions of {B}oole's calculus of logic.
\newblock {\em Mem. Am. Acad. Arts Sci.}, 9(2):317--378, 1873.
\newblock \href {https://doi.org/10.2307/25058006} {\path{doi:10.2307/25058006}}.

\bibitem{purkert}
W.~Purkert.
\newblock The double life of {F}elix {H}ausdorff/{P}aul {M}ongr\'e.
\newblock {\em Math. Intelligencer}, 30(4):36--50, 2008.
\newblock \href {https://doi.org/10.1007/BF03038095} {\path{doi:10.1007/BF03038095}}.

\bibitem{queiro}
J.~F. Queir\'o.
\newblock Jos\'e {A}nast\'acio da {C}unha: a forgotten forerunner.
\newblock {\em Math. Intelligencer}, 10(1):38--43, 1988.
\newblock \href {https://doi.org/10.1007/BF03023850} {\path{doi:10.1007/BF03023850}}.

\bibitem{quine}
W.~V.~O. Quine.
\newblock {\em Word and Object}.
\newblock MIT Press, Cambridge, new edition, 2013.
\newblock First edition published 1960.
\newblock \href {https://doi.org/10.7551/mitpress/9636.001.0001} {\path{doi:10.7551/mitpress/9636.001.0001}}.

\bibitem{reid}
C.~Reid.
\newblock {\em Courant}.
\newblock Springer, New York, 1996.
\newblock Reprint of the 1976 original.
\newblock \href {https://doi.org/10.1007/978-0-387-21626-3} {\path{doi:10.1007/978-0-387-21626-3}}.

\bibitem{rohrbach}
H.~Rohrbach.
\newblock Helmut {H}asse and {C}relle's journal.
\newblock {\em J. reine angew. Math.}, 500:5--13, 1998.
\newblock Translated from the 1982 German original by B\"arbel Deninger.
\newblock \href {https://doi.org/10.1515/crll.1998.070} {\path{doi:10.1515/crll.1998.070}}.

\bibitem{rowe}
D.~E. Rowe.
\newblock Transforming tradition: {R}ichard {C}ourant in {G}\"ottingen.
\newblock {\em Math. Intelligencer}, 37(1):20--29, 2015.
\newblock \href {https://doi.org/10.1007/s00283-014-9522-9} {\path{doi:10.1007/s00283-014-9522-9}}.

\bibitem{rudin}
W.~Rudin.
\newblock {\em Principles of Mathematical Analysis}.
\newblock McGraw-Hill, New York, first edition, 1953.

\bibitem{russ}
S.~Russ.
\newblock Bolzano's analytic programme.
\newblock {\em Math. Intelligencer}, 14(3):45--53, 1992.
\newblock \href {https://doi.org/10.1007/BF03025869} {\path{doi:10.1007/BF03025869}}.

\bibitem{russell}
B.~Russell.
\newblock Sur la logique des relations avec des applications {\`a} la th{\'e}orie des s{\'e}ries.
\newblock {\em Rev. Math.}, 7:115--148, 1901.
\newblock URL: \url{https://books.google.com/books?id=6J9FAQAAMAAJ}.

\bibitem{sierpinski}
W.~Sierpi{\'n}ski.
\newblock {\em Introduction to General Topology}.
\newblock Univ. Press, Toronto, first edition, 1934.
\newblock Translated by {C}. {Cecilia} {Krieger}.

\bibitem{simmons}
C.~K. Simmons.
\newblock Felix {H}ausdorff: ``{We} wish for you better times''.
\newblock {\em Math. Intelligencer}, 37(1):64--77, 2015.
\newblock \href {https://doi.org/10.1007/s00283-014-9474-0} {\path{doi:10.1007/s00283-014-9474-0}}.

\bibitem{sonar}
T.~Sonar.
\newblock Brunswick's second mathematical star: {R}ichard {D}edekind (1831--1916).
\newblock {\em Math. Intelligencer}, 34(2):63--67, 2012.
\newblock \href {https://doi.org/10.1007/s00283-012-9285-0} {\path{doi:10.1007/s00283-012-9285-0}}.

\bibitem{spencer}
H.~Spencer.
\newblock {\em Principles of Psychology}.
\newblock Longman, London, first edition, 1855.
\newblock URL: \url{https://oll.libertyfund.org/titles/spencer-the-principles-of-psychology-1855}.

\bibitem{vailati-reflex}
G.~Vailati.
\newblock Le propriet\`a fondamentali delle operazioni della logica deduttiva.
\newblock {\em Riv. Mat.}, I:127--134, 1891.
\newblock URL: \url{https://books.google.com/books?id=AA4MAAAAYAAJ}.

\bibitem{vailati}
G.~Vailati.
\newblock Dipendenza fra le propriet\`a delle relazioni.
\newblock {\em Riv. Mat.}, II:161--164, 1892.
\newblock URL: \url{https://books.google.com/books?id=PQ4MAAAAYAAJ}.

\bibitem{dantzig}
D.~Van~Dantzig.
\newblock Zur topologischen {A}lgebra.
\newblock {\em Math. Ann.}, 107:587--626, 1933.
\newblock \href {https://doi.org/10.1007/BF01448911} {\path{doi:10.1007/BF01448911}}.

\bibitem{whittaker}
E.~T. Whittaker.
\newblock {\em A Course of Modern Analysis}.
\newblock Univ. Press, Cambridge, first edition, 1902.
\newblock URL: \url{https://books.google.com/books?id=FSAWAAAAYAAJ}.

\bibitem{zeeman}
E.~C. Zeeman.
\newblock What's wrong with {E}uclid {B}ook {V}.
\newblock {\em Bull. Lond. Math. Soc.}, 40(1):1--17, 2008.
\newblock \href {https://doi.org/10.1112/blms/bdm104} {\path{doi:10.1112/blms/bdm104}}.

\end{thebibliography}

\end{document}